\documentclass[12pt]{article}
\usepackage{amssymb}
\usepackage{latexsym}
\usepackage[all]{xy}

\usepackage[dvips]{graphicx}
\usepackage{amssymb,color}
\usepackage[centertags]{amsmath}
\usepackage{amsfonts}

\newtheorem{thm}{Theorem}[section]

\newtheorem{lem}{Lemma}[section]

\def\BBox{\kern  -0.2cm\hbox{\vrule width 0.2cm height 0.2cm}}

\begin{document}
%
%

\title{On the Size of Equifacetted Semi-regular Polytopes}

\author{
Toma\v{z} Pisanski
\thanks{Supported in part by ARRS Grant P1Ð0294, by ESF grant EUROGiga/GReGAS, and
by a grant from the Picker Interdisciplinary Science Institute, Colgate University.}\\
Faculty of Mathematics and Physics\\
University of Ljubljana\\
Ljubljana, Slovenia\\
\and and\\[.07in]
Egon Schulte\thanks{Supported by NSF-grant DMS--0856675.}\\
Department of Mathematics\\
Northeastern University\\
Boston, Massachusetts, USA, 02115\\
\and and\\[.07in]
Asia Ivi\'{c} Weiss\thanks{Supported in part by NSERC.}\\
Department of Mathematics and Statistics\\
York University\\
Toronto, ON, Canada M3J 1P3\\}
\date{ \today }
\maketitle

\begin{abstract}
\noindent
Unlike the situation in the classical theory of convex polytopes, there is a wealth of semi-regular abstract polytopes, including interesting examples exhibiting some unexpected phenomena. We prove that even an equifacetted semi-regular abstract polytope can have an arbitrary large number of flag orbits or face orbits under its combinatorial automorphism group. 
\medskip

\end{abstract}

\section{Introduction}
\label{intro}

The study of highly-symmetric polytopes has a long history (see Coxeter~\cite{crp}). In the traditional theory, the stringent requirements included in the definition of regularity of a convex polytope can be relaxed in several different ways, yielding a great variety of weaker regularity notions (for a survey see, for example, \cite{mar}). Of particular interest are the {\em semi-regular convex polytopes\/}, which have regular facets and a vertex-transitive euclidean symmetry group. In ordinary $3$-space, the semi-regular convex polyhedra comprise the Platonic solids, two infinite classes of prisms and antiprisms, as well as the thirteen polyhedra known as {\em Archimedean solids\/} (see~\cite{cox2,cox3,john1}). In addition to the traditional {\em regular convex polytopes\/} there are just seven other semi-regular convex polytopes in higher dimensions $n$:\  three for $n = 4$, and one for each of $n = 5, 6, 7, 8$ (see \cite{blind,gos}). Arguably the most spectacular semi-regular polytopes are those in dimensions $6$, $7$ and $8$, which are related to the exceptional Coxeter groups $E_6$, $E_7$ and $E_8$ (see~\cite{crp}). A semi-regular convex polytope which is not a regular polytope has either two or three congruence (in fact, isomorphism) classes of facets. The semi-regular polytopes belong to the larger class of uniform polytopes (see \cite{clm,john1,monsch}). 

In contrast to what happens in the classical theory, there is a wealth of semi-regular abstract polytopes, including examples exhibiting some unexpected phenomena. The purpose of this note is to show that even an equifacetted semi-regular abstract polytope can have an arbitrary large number of flag orbits and face orbits under its combinatorial automorphism group. 

\section{Basic notions}
\label{}

Here we briefly introduce a few basic definitions and notions on polytopes. For a more detailed account of the theory of abstract polytopes the reader is referred to \cite{arp}. 

An (\emph{abstract\/}) \emph{polytope of rank\/} $n$, or simply an \emph{$n$-polytope\/}, is a partially ordered set $\mathcal{P}$ with a strictly monotone rank function with range $\{-1,0, \ldots, n\}$. An element of rank $j$ is a \emph{$j$-face\/} of $\mathcal{P}$, and a face of rank $0$, $1$ or $n-1$ is a \emph{vertex\/}, \emph{edge\/} or \emph{facet\/}, respectively. The maximal chains, or \emph{flags}, of $\mathcal{P}$ all contain exactly $n + 2$ faces, including a unique least face $F_{-1}$ (of rank $-1$) and a unique greatest face $F_n$ (of rank $n$).
Two flags are said to be \emph{adjacent} ($j$-\emph{adjacent}) if they differ in only one face (just their $j$-face, respectively). We shall assume that $\mathcal{P}$ is \emph{strongly flag-connected}, in the sense that, if $\Phi$ and $\Psi$ are two flags, then they can be joined by a sequence of successively adjacent flags $\Phi = \Phi_0,\Phi_1,\ldots,\Phi_k = \Psi$, each containing $\Phi \cap \Psi$. Finally, we also require that $\mathcal{P}$ has the following homogeneity property:\ whenever $F \leq G$, with $F$ a $(j-1)$-face and $G$ a $(j+1)$-face of $\mathcal{P}$ for some $j$, then there are exactly two $j$-faces $H$ of $\mathcal{P}$ with $F \leq H \leq G$. 

Recall that the \emph{order complex} of an $n$-polytope $\mathcal{P}$ is the (abstract) $(n-1)$-dimensional simplicial complex whose vertices are the proper faces of $\mathcal{P}$ (of ranks $0,\ldots,n-1$) and whose simplices are the chains (subsets of flags) which do not contain an improper face (of rank $-1$ or $n$). 

For any two faces $F$ of rank $j$ and $G$ of rank $k$ with $F \leq G$, we call
$G/F := \{ H \in \mathcal{P}\, | \, F \leq H \leq G \}$ a \emph{section} of $\mathcal{P}$, and note that $G/F$ is a polytope of rank ($k-j-1$). In particular, a face $F$ can be identified with the section $F/F_{-1}$. We also define $F_{n}/F$ to be the \emph{co-face at\/} $F$, or the \emph{vertex-figure at\/} $F$ if $F$ is a vertex.

We say that an abstract polytope $\mathcal{P}$ is {\em vertex-describable} if its faces are uniquely determined by their vertex-sets. A polytope is vertex-describable if and only if its underlying partially ordered set (of faces) can be represented by a family of subsets of the vertex-set ordered by inclusion. If a polytope $\mathcal{P}$ is a lattice, then $\mathcal{P}$ is vertex-describable. For example, the torus map $\{4,4\}_{(s,0)}$ is vertex-describable if and only if $s\geq 3$.

In this paper we shall make use of $n$-polytopes that have isomorphic facets and all vertex-figures isomorphic as well. If the facets of $\mathcal{P}$ are isomorphic to ${\cal P}_1$ and the vertex-figures are isomorphic to ${\cal P}_2$, we say that $\cal P$ is of {\em type\/} $\{{\cal P}_1,{\cal P}_2\}$ (this is a change of terminology from \cite{arp}).

Certain classes of polytopes can be described by what is known in the classical theory of polytopes as the Schl\"{a}fli symbol. We introduce it for abstract polytopes as follows. An $n$-polytope $\mathcal{P}$ is said to be {\em equivelar} if, for each $j=1,\ldots,n-1$, there exists a number $p_j$ such that, for each flag $\Psi=\{G_{-1},G_0,\ldots,G_n\}$ of $\mathcal{P}$, the section $G_{j+1}/G_{j-2}$ is a $p_{j}$-gon; in this case $\mathcal{P}$ is {\em equivelar of type\/} $\{p_1,\ldots,p_{n-1}\}$.

Highly symmetric polytopes, in particular regular polytopes, have been of interest for a long time and have also inspired many recent publications. In a sense, regular polytopes are maximally symmetric. More precisely, a polytope $\mathcal{P}$ is said to be \emph{regular\/} if its \textit{automorphism group\/} $\Gamma(\mathcal{P})$ (group of incidence preserving bijections) is transitive on flags. It easily follows that regular polytopes are equivelar.

The flag orbits of an arbitrary polytope $\mathcal{P}$ under its automorphism group $\Gamma(\mathcal{P})$ all have the same number of elements given by the order of $\Gamma(\mathcal{P})$. This follows from the fact that a polytope automorphism is uniquely determined by its effect on any flag. Thus when $\mathcal{P}$ is finite, the number of flag orbits is just the quotient of the number of flags by the order of $\Gamma(\mathcal{P})$.

An abstract polytope $\mathcal{P}$ is said to be {\em combinatorially asymmetric\/}, or simply {\em asymmetric\/}, if $\Gamma(\mathcal{P})$ is the trivial group. 

We call an abstract polytope $\mathcal{P}$ {\em equifacetted\/} if its facets are mutually isomorphic. Examples of equifacetted polytopes are given by the {\em simplicial\/} polytopes, which have facets isomorphic to simplices. All regular polytopes are equifacetted. Note that an equifacetted polytope may have a trivial automorphism group even if its facets have a large automorphism group.

In the classical theory, a semi-regular polytope is made up of regular facets and has a vertex-transitive symmetry group. This particular class of {\em geometric polytopes\/} includes the classical convex regular polytopes and {\em star-polytopes\/}, as well as the Archimedian polyhedra. We extend the classical definition in a natural way to abstract polytopes, by saying that $\mathcal{P}$ is ({\em combinatorially}) {\em semi-regular} if $\mathcal{P}$ has regular facets and $\Gamma(\mathcal{P})$ is vertex-transitive.

In this paper we only deal with equifacetted semi-regular polytopes $\mathcal{P}$. Thus the facets of $\mathcal{P}$ are mutually isomorphic regular polytopes.

\section{The $2^{\mathcal{K}}$ construction}

In this section we briefly review the polytopes $2^\mathcal{K}$ (see \cite{esext} and \cite[Section 8D]{arp}). Their discovery is originally due to Danzer; the construction announced in~\cite{dan} was never published by Danzer and first appeared in print in~\cite{esext}. The polytopes $2^\mathcal{K}$ are generalized cubes; in fact, $2^\mathcal{K}$ is just the $v$-cube when $\mathcal{K}$ is the $(v-1)$-simplex (with $v$ vertices). 

Let $\mathcal{K}$ be a finite abstract $(n-1)$-polytope with $v$ vertices and vertex-set $V:=\{1,\ldots,v\}$ (say). Suppose that $\mathcal{K}$ is vertex-describable. Then $\mathcal{P}:=2^\mathcal{K}$ will be an abstract $n$-polytope with $2^v$ vertices, each with a vertex-figure isomorphic to $\mathcal{K}$. The vertex-set of $\mathcal{P}$ is 
\begin{equation}
\label{twov}
2^V := \bigotimes_{i=1}^{v} \{0,1\} ,
\end{equation}
the cartesian product of $v$ copies of the $2$-element set $\{0,1\}$. We write elements of $2^V$ in the form $\varepsilon:=(\varepsilon_1,\ldots,\varepsilon_v)$. Now, identifying faces of $\mathcal{K}$ with their vertex-sets (recall here that $\mathcal{K}$ is vertex-describable) we take as $j$-faces of $\mathcal{P}$, for any $(j-1)$-face $F$ of $\mathcal{K}$ and any $\varepsilon$ in $2^{V}$, the subsets $F(\varepsilon)$ of $2^V$ defined by
\begin{equation}
\label{fep}
F(\varepsilon) := \{(\eta_1,\ldots,\eta_{v})\in 2^V\! \mid \eta_{i} = \varepsilon_{i} \mbox{ if } i\not\in F\}, 
\end{equation}
or, abusing notation, by the cartesian product
\[ F(\varepsilon) := \left( \bigotimes_{i \in F} \{0,1\} \right) 
\times \left( \bigotimes_{i \not\in F} \{\varepsilon_i\} \right). \]
Thus, $F(\varepsilon)$ consists of the vertices of $2^\mathcal{K}$ that coincide with $\varepsilon$ precisely in the components determined by vertices of $\mathcal{K}$ not in $F$. Then, if $F$, $F'$ are faces of $\mathcal{K}$ and $\varepsilon=(\varepsilon_1,\ldots,\varepsilon_v)$, $\varepsilon' =(\varepsilon_1,\ldots,\varepsilon_v)$ are points in $2^{V}$, we have $F(\varepsilon) \subseteq F'(\varepsilon')$ if and only if $F \leq F'$ in $\mathcal{K}$ and $\varepsilon_i = \varepsilon_{i}'$ for each $i$ which is not (a vertex) in $F'$. It is straightforward to show that the set of all faces $F(\varepsilon)$, with $F$ a face of $\mathcal{K}$ and $\varepsilon$ in $2^V$, partially ordered by inclusion (and supplemented by the empty set as least face), is an abstract $n$-polytope. This is our polytope $\mathcal{P}$. Note that the vertices $\varepsilon$ of $\mathcal{P}$ arise as $F(\varepsilon)$ with $F=\emptyset$ (that is, $F$ is the face of $\mathcal{K}$ of rank $-1$); of course, technically, $F(\varepsilon)=\{\varepsilon\}$. 

\begin{figure}[tbh]
\begin{center}
\quad
\scalebox{.50}{\includegraphics{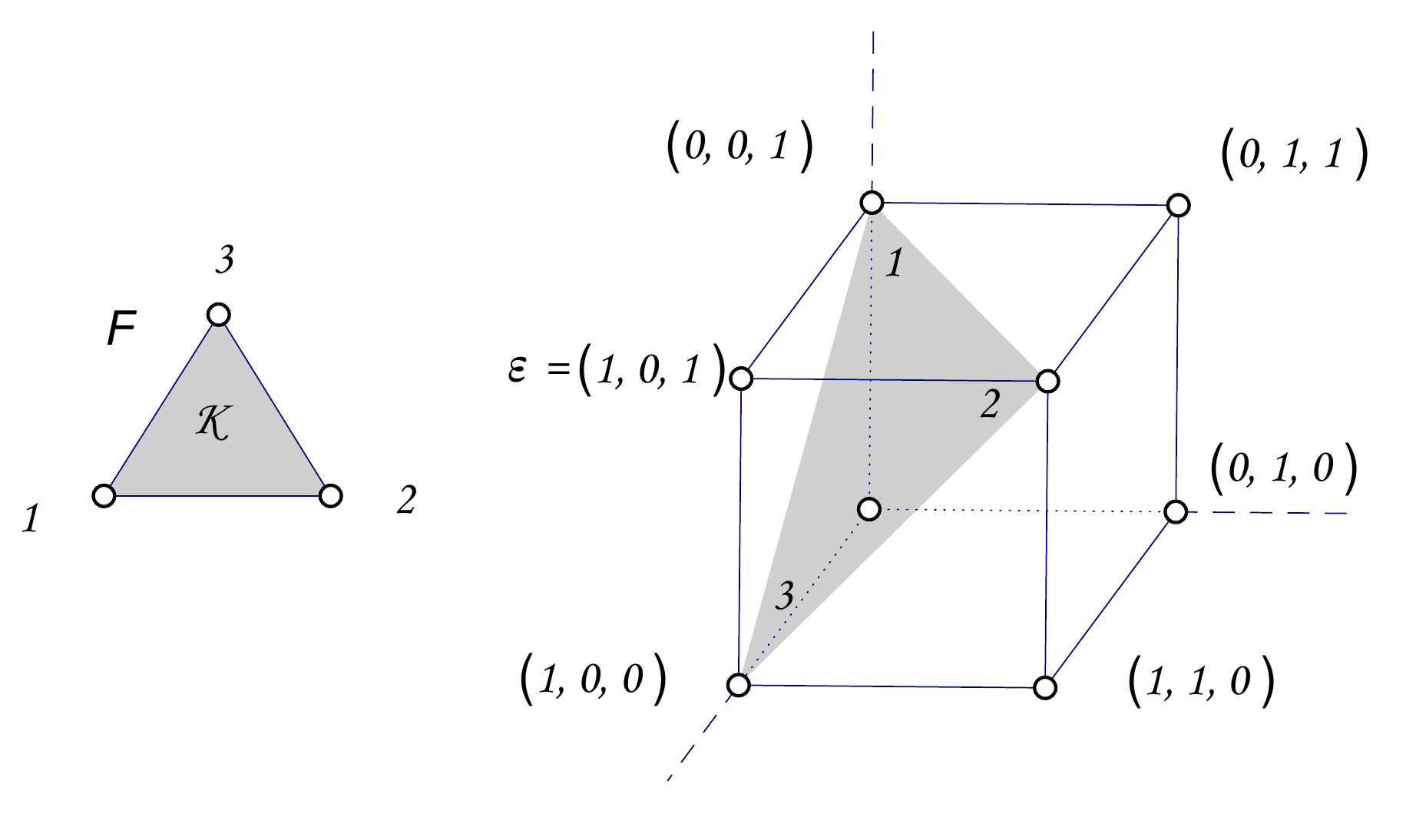}}
\caption{The $3$-cube as a polytope $2^\mathcal{K}$ obtained when $\mathcal{K}$ is a triangle. } \label{twok}
\end{center}
\end{figure}

Figure~\ref{twok} illustrates the construction of the ordinary $3$-cube as a polytope $2^\mathcal{K}$ obtained when $\mathcal{K}$ is a triangle (with vertices $1,2,3$). The triangular vertex-figure at the particular vertex $\varepsilon=(1,0,1)$ is indicated. The $2$-face $F(\varepsilon)$ determined by the edge $F=\{1,3\}$ of the triangle is the square face of the cube lying in the $xz$-plane.

The following theorem summarizes a number of important properties of $2^\mathcal{K}$.

\begin{thm}
\label{prop2k}
Let $\mathcal{K}$ be a finite abstract $(n-1)$-polytope with $v$ vertices and vertex-set~$V=\{1,\ldots,v\}$, and let $\mathcal{K}$ be vertex-describable. Then $\mathcal{P}:=2^\mathcal{K}$ has the following properties.\\[.04in]
(a) $\mathcal{P}$ is an abstract $n$-polytope with vertex-set $2^V$, and the vertex-figure at each vertex of $\mathcal{P}$ is isomorphic to $\mathcal{K}$.\\[.02in]
(b) If $F$ is a $(j-1)$-face of $\mathcal{K}$ and $\mathcal{F}:=F/F_{-1}$ (the isomorphism type of $F$ as a $(j-1)$-polytope), then each $j$-face $F(\varepsilon)$ with $\varepsilon$ in $2^V$ is isomorphic to $2^\mathcal{F}$.\\[.02in]
(c) $\Gamma(\mathcal{P}) \cong C_{2}\wr \Gamma(\mathcal{K})$, the wreath product of $C_2$ and~$\Gamma(\mathcal{K})$ defined by the natural action of $\Gamma(\mathcal{K})$ on the vertex-set of $\mathcal{P}$; in particular, $\Gamma(\mathcal{P})\cong C_{2}^{v}\rtimes \Gamma(\mathcal{K})$, a semi-direct product of the elementary abelian group $C_{2}^v$ by $\Gamma(\mathcal{K})$.\\[.02in]
(d) $\Gamma(\mathcal{P})$ acts vertex-transitively on $\mathcal{P}$, and the stabilizer of a vertex is isomorphic to $\Gamma(\mathcal{K})$.\\[.02in]
(e) If $\mathcal{K}$ is regular, then $\mathcal{P}$ is regular.
\end{thm}

\begin{proof}
For regular polytopes $\mathcal{K}$ these facts are well-known, so in particular this establishes the last part (see \cite{esext} and \cite[Section 8D]{arp}). 

For part (a), first make the following basic observation about inclusion of faces of $\mathcal{P}$ which follows immediately from the definitions:\ if $F(\varepsilon) \subseteq F'(\varepsilon')$, with $F$, $F'$, $\varepsilon$ and $\varepsilon'$ as above, then  $F'(\varepsilon') = F'(\varepsilon)$. In other words, in designating the larger face we may take $\varepsilon'=\varepsilon$.  Then, in particular, every face containing a given vertex $\varepsilon$ must necessarily be of the form $F(\varepsilon)$, with $F$ a face of $\mathcal{K}$, and any two such faces $F(\varepsilon)$ and $F'(\varepsilon)$ are incident in $\mathcal{P}$ if and only if $F$ and $F'$ are incident in $\mathcal{K}$.  This proves that the vertex-figures of $\mathcal{P}$ are isomorphic to $\mathcal{K}$. 

For part (b), if $F(\varepsilon)$ is a $j$-face of $\mathcal{P}$ and $F'(\varepsilon')$ is a face with $F'(\varepsilon')\subseteq F(\varepsilon)$, then $F'\leq F$ in $\mathcal{K}$ and $\varepsilon_{i}' = \varepsilon_{i}$ for each $i$ not in $F$; in other words, the points $\varepsilon$ and $\varepsilon'$ agree on the components representing vertices $i$ outside $F$. Hence, if we drop the components representing vertices outside $F$ and write $\eta_{F}:=(\eta_i)_{i\in F}$ for the ``trace" of a point $\eta$ on $F$, then we may safely designate the faces $F(\varepsilon)$ and $F'(\varepsilon')$ by $F(\varepsilon_F)$ and $F'(\varepsilon_{F}')$, respectively. In particular, if we now write $V_F$ for the vertex-set of $F$ and 
\[ 2^{V_F}:=\bigotimes_{i\in V_F} \{0,1\} , \] 
then $F(\varepsilon_F)= 2^{V_F}$, which is just the largest face (of rank $j$) of $2^\mathcal{F}$, and $F'(\varepsilon_{F}')$ is just a face of $2^\mathcal{F}$. Moreover, the partial order on the $j$-face $F(\varepsilon)$ of $\mathcal{P}$ is just the standard inclusion of faces in $2^\mathcal{F}$. Hence $F(\varepsilon)$ is isomorphic to $2^\mathcal{F}$.

The automorphism group $\Gamma(\mathcal{P})$ always contains an elementary abelian group $C_{2}^v$, irrespective of the symmetry properties of $\mathcal{K}$. For $k=1,\ldots,v$ let $\sigma_{k}: 2^V \rightarrow 2^V$ be the mapping that changes a point $\eta$ precisely in its $k$-th component from $0$ to $1$ or $1$ to $0$ (while leaving all other components unchanged). Clearly, $\sigma_k$ induces an automorphism of $\mathcal{P}$, which is also denoted by $\sigma_k$. Furthermore, the subgroup $\langle\sigma_1,\ldots,\sigma_v\rangle$ of $\Gamma(\mathcal{P})$ is isomorphic to $C_{2}^v$ and acts simply vertex-transitively on $\mathcal{P}$. In particular, $\Gamma(\mathcal{P})$ acts vertex-transitively.

Now consider the stabilizer of a single vertex, $\varepsilon=(0,\ldots,0)$ (say), of $\mathcal{P}$ in $\Gamma(\mathcal{P})$. Clearly, any element in this stabilizer induces an automorphism of the vertex-figure of $\mathcal{P}$ at $\varepsilon$ and hence corresponds to an automorphism of $\mathcal{K}$. Conversely, any automorphim $\varphi$ of $\mathcal{K}$ determines an automorphism $\widehat{\varphi}$ of $\mathcal{P}$ as follows. First, for an arbitrary vertex $\eta=(\eta_1,\ldots,\eta_v)$, set  
\[\widehat{\varphi}(\eta) := (\eta_{\varphi^{-1}(1)},\ldots,\eta_{\varphi^{-1}(v)}) =:\eta_\varphi, \]
and then, for an arbitrary face $F(\eta)$, define
\[ \widehat{\varphi}(F(\eta)) := \varphi(F)(\eta_\varphi) .\]
It is straightforward to check that $\widehat{\varphi}$ is indeed an automorphism of $\mathcal{P}$, which also fixes $\epsilon$. In this way, $\Gamma(\mathcal{K})$ becomes a subgroup of $\Gamma(\mathcal{P})$. In particular, $\Gamma(\mathcal{P}) \cong C_{2}\wr \Gamma(\mathcal{K}) \cong C_{2}^{v}\rtimes \Gamma(\mathcal{K})$. This establishes part (c) and also concludes the proof of part (d).
\end{proof}

\section{On the number of flag orbits}
\label{flagorb}

In this section we show that, in contrast to what happens in the classical theory, even an equifacetted semi-regular polytope can have an arbitrary large number of flag orbits under its automorphism group. 

We require the following lemma that is directly implied by the vertex-transitivity of the automorphism group.

\begin{lem}
\label{numflag}
Let $\mathcal{P}$ be any vertex-transitive $n$-polytope, and let $F$ be a vertex of $\mathcal{P}$ with vertex-figure $\mathcal{K}$. Then the number of flag orbits of $\mathcal{P}$ under $\Gamma(\mathcal{P})$ is the same as the number of flag orbits of $\mathcal{K}$ under the stabilizer $\Gamma_F(\mathcal{P})$ of $F$ in $\Gamma(\mathcal{P})$. In particular, if $\mathcal{K}$ is asymmetric, then the number of flag orbits of $\mathcal{P}$ under $\Gamma(\mathcal{P})$ is the same as the number of flags of $\mathcal{K}$.
\end{lem}

We now employ the construction of polytopes described in the previous section. In Lemma~\ref{numflag}, if $\mathcal{P}$ is a polytope of the form $2^\mathcal{K}$, then Theorem~\ref{prop2k}(d) says that the vertex-stabilizer $\Gamma_F(\mathcal{P})$ is isomorphic to $\Gamma(\mathcal{K})$. Hence, in passing from $\mathcal{K}$ to $2^\mathcal{K}$, the number of flag orbits under the automorphism group does not change.

\begin{thm}
\label{quest1}
Let $\mathcal{K}$ be an asymmetric $(n-1)$-polytope with mutually isomorphic regular facets $\mathcal{F}$, and let $\mathcal{K}$ be vertex-describable. Let $v$ denote the number of vertices and $f$ the number of flags of $\mathcal{K}$. Then $\mathcal{P}:= 2^{\mathcal{K}}$ is an equifacetted semi-regular $n$-polytope with facets isomorphic to $2^{\mathcal{F}}$, automorphism group $\Gamma(\mathcal{P}) = C_{2}^{v}$, and $f$ flag orbits under $\Gamma(P)$.
\end{thm}

\begin{proof}
Recall from Theorem~\ref{prop2k} of the previous section that $\mathcal{P}$ is vertex-transitive and that, since $\mathcal{F}$ is regular, $\mathcal{P}$ has regular facets isomorphic to $2^{\mathcal{F}}$. Hence $\mathcal{P}$ is an equifacetted semi-regular $n$-polytope. Moreover, since $\Gamma(\mathcal{K})$ is trivial, we have $\Gamma(\mathcal{P}) = C_{2}^{v}\rtimes \Gamma(\mathcal{K}) = C_{2}^{v}$.  By the previous lemma, $f$ is the number of flag orbits of $\mathcal{P}$ under $\Gamma(\mathcal{P})$.
\end{proof}

Note that the condition that $\mathcal{K}$ be vertex-describable is rather weak and is satisfied for most polytopes. 

\begin{figure}[tbh]
\begin{center}
\quad
\scalebox{.50}{\includegraphics{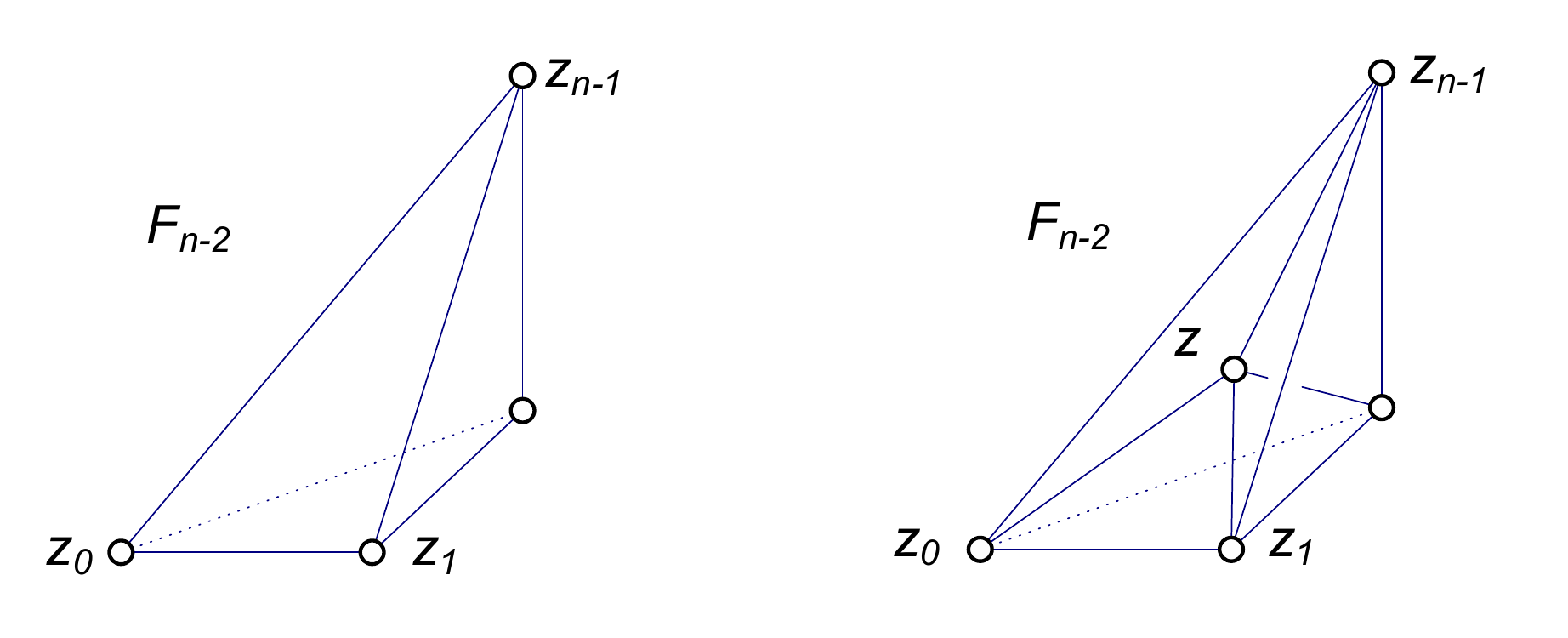}}
\caption{Subdivision of a facet of $\mathcal{L}'$} \label{pointz}
\end{center}
\end{figure}

There are numerous polytopes $\mathcal{K}$ satisfying the assumptions of the above theorem, including many simplicial polytopes (in fact, even many simplicial convex polytopes). Examples (of abstract polytopes) can be obtained by the following construction, which also demonstrates an interesting method of local symmetry-breaking. Starting with an equivelar $(n-1)$-polytope $\mathcal{L}$ of type $\{p_1,\ldots,p_{n-2}\}$ with $p_j\geq 3$ for each $j$, let $\mathcal{L}'$ denote its order complex. Assume that $p_1,\ldots, p_{n-2}$ are mutually distinct.  In $\mathcal{L}'$, replace a single simplicial facet, $F_{n-2}$ (say), by $n-1$ simplices all sharing a ``central" vertex $z$ (say) not contained in $\mathcal{L}'$, as indicated in Figure~\ref{pointz}, while keeping all other faces of $\mathcal{L}'$ invariant. (Recall here that $F_{n-2}$ is a flag of the original polytope $\mathcal{L}$.) The resulting polytope $\mathcal{K}$ is again simplicial and has just one $(n-1)$-valent vertex, namely $z$.  We will show that $\mathcal{K}$ is asymmetric. 

Let $\varphi$ be any automorphism of $\mathcal{K}$. Then $\varphi$ must necessarily fix $z$ and belong to the stabilizer of $F_{n-2}$ in $\Gamma(\mathcal{L}')$. Denote the vertices of $F_{n-2}$ by $z_0,\ldots,z_{n-1}$, where $z_j$ is the vertex corresponding to the $j$-face of the original polytope $\mathcal{L}$ in the flag corresponding to $F_{n-2}$. In $\mathcal{L}'$, the $(n-3)$-face $G_{i,j}$ of $F_{n-2}$ with vertex set $\{z_0,\ldots,z_{n-1}\}\setminus \{z_i,z_j\}$ is surrounded by $p_{ij}$ facets, where $p_{ij}=4$ if $|i-j|\geq 2$, and $p_{ij}=p_{i+1}$ if $j=i+1$. Since the $p_{ij}$ are distinct, each $(n-3)$-face $G_{i,i+1}$ of $F_{n-2}$ is completely determined by the number of facets of $\mathcal{L}$ that surround it. Hence $\varphi$ must fix every face $G_{i,i+1}$. It follows that $\varphi$ must fix every vertex of $F_{n-1}$. Thus, by the connectedness properties, $\varphi$ must be the identity isomorphism of $\mathcal{L}'$ and hence of $\mathcal{K}$.

When $n-1=3$ we may take $\mathcal{L}$ to be any toroidal map of type $\{3,6\}$ and construct a $3$-polytope $\mathcal{K}$ with triangular (regular) facets and trivial automorphism group. Note that, by our choice of the maps, the number of flags $f$ of $\mathcal{K}$ can be taken to be arbitrarily large. Inductively, when $n-1\geq 4$ we can appeal to an extension theorem for regular polytopes, proved in Pellicer~\cite{pel}, which establishes the existence of finite regular polytopes with arbitrarily preassigned isomorphism type of facets and an arbitrarily preassigned even number for the last entry of the Schl\"afli symbol. Hence, starting from regular maps of type $\{3,6\}$, we find infinitely many regular polytopes $\mathcal{L}$, each with distinct entries in the 
Schl\"afli symbol, which then, via the order complex $\mathcal{L}'$, provide polytopes $\mathcal{K}$ with the desired properties. 

In conclusion, we have established that there are equifacetted semi-regular polytopes with an arbitrary large number of flag orbits under their automorphism group. 

\section{On the number of face orbits}
\label{facetorb}

In this section we establish that an equifacetted semi-regular polytope can also have an arbitrary large number of $j$-face orbits under its automorphism group, for each $j=0,\ldots,n-1$. A priori, this is not implied by the results of the previous section. 

We begin with the following lemma for face orbits, which is weaker than the corresponding Lemma~\ref{numflag} for flag orbits. 

\begin{lem}
\label{numfacet}
Let $\mathcal{P}$ be any vertex-transitive $n$-polytope, let $F$ be a vertex of $\mathcal{P}$ with vertex-figure $\mathcal{K}$, and let $1\leq j\leq n-1$. Then the number of $j$-face orbits of $\mathcal{P}$ under $\Gamma(\mathcal{P})$ is at most the number of $(j-1)$-face orbits of $\mathcal{K}$ under the stabilizer $\Gamma_F(\mathcal{P})$ of $F$ in $\Gamma(\mathcal{P})$. In particular, if $\mathcal{K}$ is asymmetric, then the number of $j$-face orbits of $\mathcal{P}$ under $\Gamma(\mathcal{P})$ is at most the number of $(j-1)$-faces of $\mathcal{K}$.
\end{lem}

\begin{proof}
By the vertex-transitivity of $\mathcal{P}$, each $j$-face orbit of $\mathcal{P}$ contains a $j$-face of which $F$ is a vertex. For a $j$-face $G$ containing $F$, let $G^{\Gamma(\mathcal{P})}$ and $G^{\Gamma_{F}(\mathcal{P})}$ denote the orbit of $G$ under $\Gamma(\mathcal{P})$ or $\Gamma_{F}(\mathcal{P})$ respectively. Then $G^{\Gamma_{F}(\mathcal{P})}\rightarrow G^{\Gamma(\mathcal{P})}$ determines a well-defined and surjective mapping $\gamma$ from the set of all $j$-face orbits of $j$-faces with vertex $F$ under $\Gamma_{F}(\mathcal{P})$ to the set of all $j$-face orbits of $\mathcal{P}$ under $\Gamma_{F}(\mathcal{P})$. This simply follows from the fact that $\Gamma_{F}(\mathcal{P})$ is a subgroup of $\Gamma(\mathcal{P})$, and proves the first part of the lemma. (Bear in mind here that the $(j-1)$-faces of $\mathcal{K}$ are just the $j$-faces of $\mathcal{P}$ containing $F$.)  For the second part, note that $\Gamma_F(\mathcal{P})$ is also trivial if $\Gamma(\mathcal{K})$ is trivial. 
\end{proof}

The next lemma establishes that for certain kinds of polytopes the inequality of the previous lemma becomes an equality.

\begin{lem}
\label{inj}
Let $\mathcal{K}$ be an asymmetric vertex-describable $(n-1)$-polytope, let $\mathcal{P}:=2^{\mathcal{K}}$, let $F$ be a vertex of $\mathcal{P}$, and let $1\leq j\leq n-1$. Then the number of $j$-face orbits of $\mathcal{P}$ under $\Gamma(\mathcal{P})$ is equal to the number of $(j-1)$-faces of $\mathcal{K}$.
\end{lem}

\begin{proof}
We only need to prove that the mapping $\gamma$ defined in the previous proof is also injective. Suppose that 
$G_{1}^{\Gamma(\mathcal{P})}=G_{2}^{\Gamma(\mathcal{P})}$, where $G_1$ and $G_2$ are $j$-faces of $\mathcal{P}$ containing the vertex $F$. Now $F=\varepsilon$ for some $\varepsilon$ in the vertex-set of $\mathcal{P}$, and the $j$-faces $G_i$ are of the form $K_{i}(\varepsilon)$ with the same $\varepsilon$ and with $K_{i}$ the vertex set of a $(j-1)$-face of $\mathcal{K}$. Since $\Gamma(\mathcal{K})$ is trivial, $\Gamma(\mathcal{P})=C_{2}^{v}$ and hence $\Gamma(\mathcal{P})$ acts only on the $\varepsilon$-component in the definition of a face. It follows that $G_1=G_2$ and hence $G_{1}^{\Gamma_{F}(\mathcal{P})}=G_{2}^{\Gamma_{F}(\mathcal{P})}$. 
\end{proof}

The next theorem now follows immediately and can be proved in the same manner as Theorem~\ref{quest1}.

\begin{thm}
\label{quest2}
Let $\mathcal{K}$ be an asymmetric vertex-describable $(n-1)$-polytope with mutually isomorphic regular facets $\mathcal{F}$, and let $1\leq j\leq n-1$. Let $v$ denote the number of vertices and $f_{j-1}$ the number of $(j-1)$-faces of $\mathcal{K}$. Then $\mathcal{P}:= 2^{\mathcal{K}}$ is an equifacetted semi-regular $n$-polytope with facets isomorphic to $2^{\mathcal{F}}$, automorphism group $\Gamma(\mathcal{P}) = C_{2}^{v}$, and $f_{j-1}$\ $j$-face orbits under $\Gamma(P)$.
\end{thm}

Now the polytopes constructed in the previous section can also be used to show that equifacetted semi-regular polytopes can have an arbitrarily large number of $j$-face orbits under their automorphism group for each $j=1,\ldots,n-1$ (in fact, even simultaneously for all these $j$). 
\vskip.1in
\noindent
{\bf Acknowledgment}
We are grateful to an anonymous referee for a number of helpful comments.

\end{document}